\documentclass[11pt]{article}

\usepackage{amsmath}
\usepackage{amsthm}
\usepackage{amssymb}
\usepackage{amscd}
\theoremstyle{definition}
\newtheorem{theorem}{Theorem}[section]
\newtheorem{prop}[theorem]{Proposition}
\newtheorem{lemma}[theorem]{Lemma}
\newtheorem{corollary}[theorem]{Corollary}
\newtheorem{definition}[theorem]{Definition}

\newtheorem{remark}[theorem]{Remark}

\numberwithin{equation}{section}
\newenvironment{demo}[1]{%
  \trivlist
  \item[\hskip\labelsep
        {\it #1.}]
}{%
\hfill\qedsymbol
  \endtrivlist
}

\newcommand\Pf{\operatorname{Pf}}
\newcommand\sgn{\operatorname{sgn}}
\newcommand\trans{{}^t\!}

\newcommand\ep{\varepsilon}

\renewcommand\tilde{\widetilde}
\renewcommand\hat{\widehat}

\allowdisplaybreaks

\title{
Generalized Sylvester Formulas and skew Giambelli Identities
}

\author{
Soichi OKADA%
\footnote{
Graduate School of Mathematics, Nagoya University, 
Furo-cho, Chikusa-ku, Nagoya 464-8602, Japan, 
{\tt okada@math.nagoya-u.ac.jp}
}
\footnote{
This work was partially supported by 
JSPS Grants-in-Aid for Scientific Research No.~24340003 and No.~15K13425.
}
}

\date{
}

\begin{document}

\maketitle

\begin{abstract}
We obtain a common generalization of two types of Sylvester formulas for compound determinants 
and its Pfaffian analogue.
As applications, we give generalizations of the Giambelli identity to skew Schur functions 
and the Schur identity to Schur's skew $Q$-functions.
\end{abstract}

%%%%%%%%%%
% Introduction
%%%%%%%%%%

\section{%
Introduction
}

Determinant and Pfaffian formulas for general matrices 
are a powerful tool in proving relations among special functions 
and in evaluating specific determinants and Pfaffians.
For example, in \cite{ItO}, we applied the Cauchy--Sylvester formula for compound determinants 
to give a transparent proof of the evaluation of the determinant involving 
$BC_n$-type Jackson integrals.
And, in \cite{IIO}, we found a variant of the Cauchy--Sylvester formula 
and obtained product evaluations of determinants of classical group characters.
The aims of this paper are to establish a common generalization of two types of Sylvester formulas 
for compound determinants and its Pfaffian analogue, and to obtain symmetric function identities 
by applying these generalized formulas.

Given a matrix $X = \bigl( x_{i,j} \bigr)$ and sequences $I = (i_1, \dots, i_p)$ of row indices 
and $I' = (i'_1, \dots, i'_q)$ of column indices, let $X \begin{pmatrix} I \\ I' \end{pmatrix}$ 
be the $p \times q$ matrix 
obtained from $X$ by picking up the rows indexed by $I$ and the columns indexed by $I'$, i.e.,
$$
X \begin{pmatrix} I \\ I' \end{pmatrix}
 =
\left( x_{i_\alpha, i'_\beta} \right)_{1 \le \alpha \le p, \, 1 \le \beta \le q}.
$$
For two sequences $I$ and $J$, we denote by $I \sqcup J$ the concatenation of $I$ and $J$, 
and, if $J$ is a subsequence of $I$, then we denote by $I \setminus J$ 
the sequence obtained by removing the entries of $J$ from $I$.
For example, $(1,3,4) \sqcup (2,6) = (1,3,4,2,6)$ and $(1,3,4,2,6) \setminus (3,2) = (1,4,6)$.
Then the Sylvester formula and its dual version are stated as follows:

\begin{prop}
\label{prop:Sylvester}
Let $X$ be a matrix and let $K$ and $K'$ be sequences of row and column indices of the same length.
\begin{enumerate}
\item[(1)]
(Sylvester \cite{Sylvester})
For sequences $I = (i_1, \dots, i_p)$ of row indices and $I' = (i'_1, \dots, i'_p)$ of column indices, 
we have
\begin{equation}
\label{eq:Sylvester1}
\det \left(
 \det X \begin{pmatrix} ( i_\alpha ) \sqcup K \\ ( i'_\beta ) \sqcup K' \end{pmatrix}
 \right)_{1 \le \alpha, \, \beta \le p}
 =
\det X \begin{pmatrix} I \sqcup K \\ I' \sqcup K' \end{pmatrix}
\cdot
\left( \det X \begin{pmatrix} K \\ K' \end{pmatrix} \right)^{p-1}.
\end{equation}
\item[(2)]
For sequences $J = (j_1, \dots, j_q)$ of row indices and $J' = (j'_1, \dots, j'_q)$ of column indices, 
we have
\begin{equation}
\label{eq:Sylvester2}
\det \left(
 \det X \begin{pmatrix} (J \setminus (j_\alpha)) \sqcup K \\ (J \setminus (j'_\beta)) \sqcup K' \end{pmatrix}
\right)_{1 \le \alpha, \, \beta \le q}
\\
 =
\left( \det X \begin{pmatrix} J \sqcup K \\ J' \sqcup K' \end{pmatrix} \right)^{q-1}
\cdot
\det X \begin{pmatrix} K \\ K' \end{pmatrix}.
\end{equation}
\end{enumerate}
\end{prop}

The identities (\ref{eq:Sylvester1}) and (\ref{eq:Sylvester2}) are dual to each other 
in the sense that one is obtained from the other by applying to the matrix $X^{-1}$ 
and using Jacobi's complementary minor formula \cite{Jacobi}. 
One of the main results of this paper is the following theorem, 
which is a common generalization of (\ref{eq:Sylvester1}) and (\ref{eq:Sylvester2}).
In fact, we can recover (\ref{eq:Sylvester1}) (resp. (\ref{eq:Sylvester2})) by specializing $q=0$ (resp. $p=0$) 
in (\ref{eq:main1}).

\begin{theorem}
\label{thm:main1}
Let $X$ be a matrix and $I$, $J$ and $K$ (resp. $I'$, $J'$ and $K'$) sequences of row (resp. column) indices 
of length $p$, $q$ and $r$.
We define a $(p+q) \times (p+q)$ matrix 
$\tilde{X} = \bigl( \tilde{x}_{\alpha, \beta} \bigr)_{1 \le \alpha, \, \beta \le p+q}$ by
$$
\begin{cases}
\tilde{x}_{\alpha, \beta}
 =
\det X \begin{pmatrix} (i_\alpha) \sqcup J \sqcup K \\ (i'_\beta) \sqcup J' \sqcup K' \end{pmatrix} 
 &\text{if $1 \le \alpha, \beta \le p$,}
\\
\tilde{x}_{\alpha, p+\beta}
 =
\det X \begin{pmatrix} (i_\alpha) \sqcup ( J \setminus (j_\beta) ) \sqcup K \\ J' \sqcup K' \end{pmatrix}
 &\text{if $1 \le \alpha \le p$ and $1 \le \beta \le q$,}
\\
\tilde{x}_{p+\alpha, \beta}
 =
\det X \begin{pmatrix} J \sqcup K \\ (i'_\beta) \sqcup ( J' \setminus (j'_\alpha) ) \sqcup K' \end{pmatrix}
 &\text{if $1 \le \alpha \le q$ and $1 \le \beta \le p$,}
\\
\tilde{x}_{p+\alpha, p+\beta}
 =
- \det X \begin{pmatrix} ( J \setminus (j_\beta) ) \sqcup K \\ (J' \setminus (j'_\alpha) ) \sqcup K' \end{pmatrix}
 &\text{if $1 \le \alpha, \beta \le q$.}
\end{cases}
$$
Then we have
\begin{equation}
\label{eq:main1}
\det \tilde{X}
 =
(-1)^q
\det X \begin{pmatrix} I \sqcup K \\ I' \sqcup K' \end{pmatrix}
\cdot
\left( \det X \begin{pmatrix} J \sqcup K \\ J' \sqcup K' \end{pmatrix} \right)^{p+q-1}.
\end{equation}
\end{theorem}

For example, if $I = I' = (1,2)$, $J = J' = (3,4)$, $K = K' = (5)$, then Equation (\ref{eq:main1}) reads
$$
\det \begin{pmatrix}
 \begin{vmatrix} 1,3,4,5 \\ 1,3,4,5 \end{vmatrix}
&
 \begin{vmatrix} 1,3,4,5 \\ 2,3,4,5 \end{vmatrix}
&
 \begin{vmatrix} 1,4,5 \\ 3,4,5 \end{vmatrix}
&
 \begin{vmatrix} 1,3,5 \\ 3,4,5 \end{vmatrix}
\\[13pt]
 \begin{vmatrix} 2,3,4,5 \\ 1,3,4,5 \end{vmatrix}
&
 \begin{vmatrix} 2,3,4,5 \\ 2,3,4,5 \end{vmatrix}
&
 \begin{vmatrix} 2,4,5 \\ 3,4,5 \end{vmatrix}
&
 \begin{vmatrix} 2,4,5 \\ 3,4,5 \end{vmatrix}
\\[13pt]
 \begin{vmatrix} 3,4,5 \\ 1,4,5 \end{vmatrix}
&
 \begin{vmatrix} 3,4,5 \\ 2,4,5 \end{vmatrix}
&
 - \begin{vmatrix} 4,5 \\ 4,5 \end{vmatrix}
&
 - \begin{vmatrix} 3,5 \\ 4,5 \end{vmatrix}
\\[13pt]
 \begin{vmatrix} 3,4,5 \\ 1,3,5 \end{vmatrix}
&
 \begin{vmatrix} 3,4,5 \\ 2,3,5 \end{vmatrix}
&
 - \begin{vmatrix} 4,5 \\ 4,5 \end{vmatrix}
&
 - \begin{vmatrix} 3,5 \\ 3,5 \end{vmatrix}
\end{pmatrix}
 =
(-1)^2
\begin{vmatrix} 1,2,5 \\ 1,2,5 \end{vmatrix}
\cdot
\begin{vmatrix} 3,4,5 \\ 3,4,5 \end{vmatrix}^{2+2-1},
$$
where $\begin{vmatrix} i_1, \dots, i_p \\ i'_1, \dots, i'_p \end{vmatrix}$ stands for 
$\det X \begin{pmatrix} i_1, \dots, i_p \\ i'_1, \dots, i'_p \end{pmatrix}$.

Knuth \cite{Knuth} gave Pfaffian analogues of the Sylvester formulas in Theorem~\ref{prop:Sylvester}.
For a skew-symmetric matrix $Y$ and a sequence $I = (i_1, \cdots, i_l)$ of row/column indices, 
we write $Y(I)$ for $Y \begin{pmatrix} I \\ I \end{pmatrix}$.

\begin{prop}
\label{prop:Knuth}
Let $Y$ be a skew-symmetric matrix and $K$ be a sequence of row/column indices of even length.
\begin{enumerate}
\item[(1)]
(Knuth \cite[2.5]{Knuth})
For a sequence $I = (i_1, \dots, i_l)$ of row/column indices of even length $l$, we have
\begin{equation}
\label{eq:Knuth1}
\Pf
 \Bigl( 
  \Pf Y( ( i_\alpha, i_\beta ) \sqcup K )
 \Bigl)_{1 \le \alpha < \beta \le l}
 =
\Pf Y(I \sqcup K)
\cdot
\bigl( \Pf Y( K ) \bigr)^{l/2-1}.
\end{equation}
\item[(2)]
(Knuth \cite[(2.7)]{Knuth})
For a sequence $J = (j_1, \dots, j_m)$ of row/column indices of even length $m$, we have
\begin{equation}
\label{eq:Knuth2}
\Pf
 \Bigl( 
  \Pf Y( (J \setminus ( j_\alpha, j_\beta ) \sqcup K )
 \Bigr)_{1 \le \alpha < \beta \le m}
 =
\bigl( \Pf Y( J \sqcup K ) \bigr)^{m/2-1}
\cdot
\Pf Y( K ).
\end{equation}
\end{enumerate}
\end{prop}

Another main result is the following common generalization of (\ref{eq:Knuth1}) and (\ref{eq:Knuth2}).
If $m=0$ (resp. $l=0$), then Equation (\ref{eq:main2}) reduces to (\ref{eq:Knuth1}) (resp. (\ref{eq:Knuth2})).
Note that Hirota \cite{Hirota} gave the special case $l=m$.

\begin{theorem}
\label{thm:main2}
Let $Y$ be a skew-symmetric matrix.
Let $l$, $m$ and $n$ be nonnegative integers with the same parity,
and let $I$, $J$ and $K$ be sequences of row/column indices of $Y$ of length $l$, $m$ and $n$ respectively.
We define an $(l+m) \times (l+m)$ skew-symmetric matrix $\tilde{Y}
 = \bigl( \tilde{y}_{\alpha, \beta} \bigr)_{1 \le \alpha, \beta \le p+q}$ by
$$
\begin{cases}
\tilde{y}_{\alpha, \beta}
 =
\Pf Y( (i_\alpha, i_\beta) \sqcup J \sqcup K )
 &\text{if $1 \le \alpha < \beta \le l$,}
\\
\tilde{y}_{\alpha, m+\beta}
 =
\Pf Y( (i_\alpha) \sqcup ( J \setminus (j_\beta)) \sqcup K )
 &\text{if $ 1 \le \alpha \le l$ and $1 \le \beta \le m$,}
\\
\tilde{y}_{m+\alpha,m+\beta}
 =
\Pf Y( (J \setminus (j_\alpha, j_\beta)) \sqcup K )
 &\text{if $1 \le \alpha < \beta \le m$.}
\end{cases}
$$
Then we have
\begin{equation}
\label{eq:main2}
\Pf \tilde{Y}
 =
\Pf Y(I \sqcup K)
\cdot
\bigl( \Pf Y(J \sqcup K) \bigr)^{(l+m)/2-1}.
\end{equation}
\end{theorem}

As an application of our generalized Sylvester formula (Theorem~\ref{thm:main1}), 
we give an extension of the Giambelli identity to skew Schur functions, 
which expresses a skew Schur function $s_{\lambda/\mu}$ 
in terms of a determinants of skew Schur functions of the form $\pm s_{\nu/\mu}$.
Similarly, we use the generalized Knuth formula (Theorem~\ref{thm:main2}) to obtain 
a generalization of the Schur identity to Schur's skew $Q$-functions,
which expresses a skew $Q$-function $Q_{\lambda/\mu}$ 
in terms of a Pfaffian of skew $Q$-functions of the form $\pm Q_{\nu/\mu}$.

This article is organized as follows.
In Section~2, we prove Theorem~\ref{thm:main2} by using a Pfaffian analogue of 
the Desnanot--Jacobi identity.
In Section~3, we derive Theorem~\ref{thm:main1} from Theorem~\ref{thm:main2} by 
using a relation between determinants and Pfaffians.
Applications to symmetric function identities are provided in Sections~4 and 5.

%%%%%%%%%%
% Proof of generalized Sylvester identity for Pfaffians
%%%%%%%%%%

\section{%
Proof of Theorem~\ref{thm:main2}
}

Recall a definition and some properties of Pfaffians (see \cite{IsO} for details).
Given an $2m \times 2m$ skew-symmetric matrix $A = \bigl( a_{ij} \bigr)_{1 \le i, j \le 2m}$, 
the Pfaffian of $A$, denoted by $\Pf A$, is defined by
$$
\Pf A
 =
\sum_{\pi \in F_{2m}}
 \sgn(\pi) a_{\pi(1), \pi(2)} a_{\pi(3), \pi(4)} \dots a_{\pi(2m-1), \pi(2m)},
$$
where $F_{2m}$ is the subset of the symmetric group $S_{2m}$ given by
$$
F_{2m}
 =
\{ \pi \in S_{2m} : \pi(1) < \pi(3) < \dots < \pi(2m-1), \, \pi(2i-1) < \pi(2i) \, (1 \le i \le m) \}.
$$
For example, if $2m=4$, then we have
$$
\Pf A = a_{12} a_{34} - a_{13} a_{24} + a_{14} a_{23}.
$$
For a $2m \times 2m$ skew-symmetric matrix $A$ and a $2m \times 2m$ matrix $T$, we have
$$
\Pf \left( \trans T A T \right) = \det T \cdot \Pf A.
$$
Hence the Pfaffian is alternating in the sense that 
$\Pf \left( a_{\sigma(i), \sigma(j)} \right) 
 = \sgn(\sigma) \Pf \left( a_{i,j} \right)$ for any permutation $\sigma \in S_{2m}$.

One key ingredient of the proof of Theorem~\ref{thm:main2} is the Pfaffian analogue of the Desnanot--Jacobi formula, 
which is a special case of the basic identity in \cite{Knuth}.

\begin{lemma}
\label{lem:Pf-Dodgeson}
(See e.g. \cite[(1.1)]{Knuth}.)
Let $A$ be a skew-symmetric matrix $A$ with rows/columns indexed by $( 1, 2, \dots, 2m )$.
For row/column indices $i < j < k < l$, we have
\begin{equation}
\label{eq:Pf-dodgeson}
\Pf A^{i,j} \cdot \Pf A^{k,l} - \Pf A^{i,k} \cdot \Pf A^{j,l} + \Pf A^{i,l} \cdot \Pf A^{j,k}
 =
\Pf A \cdot \Pf A^{i,j,k,l},
\end{equation}
where $A^{i_1, \cdots, i_p}$ denotes the skew-symmetric matrix obtained from $A$ 
by removing rows and columns with indices $i_1, \cdots, i_p$.
\end{lemma}

First we show the following identities, which are special cases of Theorem~\ref{thm:main2}.

\begin{lemma}
\label{lem:special}
Let $Y$ be a skew-symmetric matrix.
Let $(a, b, c, d)$ and $K$ be sequences of row/column indices.
\begin{enumerate}
\item[(1)]
If $k$ is even, then we have
\begin{multline}
\label{eq:special22}
\Pf Y( (a,b,c,d) \sqcup K ) \cdot \Pf (K)
 - \Pf Y( (a,d) \sqcup K ) \cdot \Pf Y( (b,c) \sqcup K )
\\
 + \Pf Y( (a,c) \sqcup K ) \cdot \Pf Y( (b,d) \sqcup K )
 =
\Pf Y( (a,b) \sqcup K ) \cdot \Pf Y( (c,d) \sqcup K ).
\end{multline}
\item[(2)]
If $k$ is odd, then we have
\begin{multline}
\label{eq:special13}
\Pf Y( ( a,c,d) \sqcup K ) \cdot \Pf Y( (b) \sqcup K)
 - \Pf Y( ( a,b,d) \sqcup K ) \cdot \Pf Y( (c) \sqcup K)
\\
 + \Pf Y( ( a,b,c) \sqcup K ) \cdot \Pf Y( (d) \sqcup K)
 =
\Pf Y( ( b,c,d) \sqcup K ) \cdot \Pf Y( (a) \sqcup K).
\end{multline}
\end{enumerate}
\end{lemma}

\begin{demo}{Proof}
The first identity (\ref{eq:special22}) is a restatement of Lemma~\ref{lem:Pf-Dodgeson}.
We shall prove the second identity (\ref{eq:special13}).
Let $e$ be the first entry of $K$ and let $K' = K \setminus (e)$ be the remaining sequence.
If we write $[i_1, \dots, i_r] = \Pf Y( (i_1, \dots, i_r) \sqcup K' )$, 
then the left hand side of (\ref{eq:special13}) is written as
$$
[a,c,d,e] \cdot [b,e] - [a,b,d,e] \cdot [c,e] + [a,b,c] \cdot [d,e].
$$
By using lemma~\ref{lem:Pf-Dodgeson}, we have
\begin{align*}
[a,c,d,e] \cdot [\emptyset] &= [a,c] \cdot [d,e] - [a,d] \cdot [c,e] + [a,e] \cdot [c,d],
\\
[a,b,d,e] \cdot [\emptyset] &= [a,b] \cdot [d,e] - [a,d] \cdot [b,e] + [a,e] \cdot [b,d],
\\
[a,b,c,e] \cdot [\emptyset] &= [a,b] \cdot [c,e] - [a,c] \cdot [b,e] + [a,e] \cdot [b,c],
\end{align*}
where $[\emptyset] = \Pf Y(K')$.
Hence the left hand side of (\ref{eq:special13}) is equal to
$$
\frac{ 1 }{ [\emptyset] }
\cdot [a,e]
\left(
 [b,c] \cdot [d,e] - [b,d] \cdot [c,e] + [b,e] \cdot [c,d]
\right),
$$
which turns out to be equal to $[a,e] \cdot [b,c,d,e] = \Pf Y ( (a) \sqcup K ) \cdot \Pf Y( (b,c,d) \sqcup K )$ 
again by using Lemma~\ref{lem:Pf-Dodgeson}.
\end{demo}

\begin{demo}{Proof of Theorem~\ref{thm:main2}}
For three sequences $I$, $J$ and $K$, we denote by $Z(I,J,K)$ the $(l+m) \times (l+m)$ skew-symmetric 
matrix given by
$$
Z(I,J,K)
=
\begin{pmatrix}
 \Bigl( 
  \Pf Y( ( i_\alpha, i_\beta ) \sqcup J \sqcup K )
 \Bigr)_{1 \le \alpha < \beta \le l}
&
 \Bigl(
  \Pf Y( ( i_\alpha ) \sqcup ( J \setminus ( j_\beta ) \sqcup K )
 \Bigr)_{\substack{1 \le \alpha \le l \\ 1 \le \beta \le m}}
\\
 &
 \Bigl(
  \Pf Y( (J \setminus ( j_\alpha, j_\beta )) \sqcup K )
 \Bigr)_{1 \le \alpha < \beta \le m}
\end{pmatrix},
$$
where $l$ and $m$ are the lengths of $I$ and $J$ respectively.
We proceed by induction on $l$ and $m$.

First we consider the case where $l \ge 4$.
In this case, we apply Lemma~\ref{lem:Pf-Dodgeson} to the matrix $A = \tilde{Y} = Z(I,J,K)$ 
with $(i,j,k,l) = (1,2,3,4)$.
By using the induction hypothesis on $l$, we have
\begin{align*}
\tilde{Y}^{\alpha, \beta}
 &=
\Pf Z(I \setminus ( i_\alpha, i_\beta ), J, K)
\\
 &=
\Pf Y( (I \setminus ( i_\alpha, i_\beta )) \sqcup K )
\cdot
\bigl ( \Pf Y( J \sqcup K ) \bigr)^{((l-2)+m)/2-1}
\quad(1 \le \alpha < \beta \le 4),
\end{align*}
and
\begin{align*}
\tilde{Y}^{1,2,3,4}
 &=
\Pf Z(I \setminus ( i_1, i_2, i_3, i_4 ), J, K)
\\
 &=
\Pf Y( (I \setminus ( i_1, i_2, i_3, i_4 ) \sqcup K )
\cdot
\bigl( \Pf Y( J \sqcup K ) \bigr)^{((l-4)+m)/2-1}.
\end{align*}
Hence we see that
\begin{align*}
&
\Pf Z(I,J,K)
\\
 &\quad
=
\frac{ 1 }
     { \Pf Y(I \setminus (i_1, i_2, i_3, i_4)) }
\left\{
 \begin{array}{l}
  \Pf Y( (I \setminus ( i_1, i_2 )) \sqcup K )
  \cdot
  \Pf Y( (I \setminus ( i_3, i_4 )) \sqcup K )
 \\
  -
  \Pf Y( (I \setminus ( i_1, i_3 )) \sqcup K )
  \cdot
  \Pf Y( (I \setminus ( i_2, i_4 )) \sqcup K )
 \\
 +
  \Pf Y( (I \setminus ( i_1, i_4 )) \sqcup K )
  \cdot
  \Pf Y( (I \setminus ( i_2, i_3 )) \sqcup K )
 \end{array}
\right\}
\\
&\quad\quad
\times
\bigl( \Pf Y(J \sqcup K) \bigr)^{(l+m)/2-1}.
\end{align*}
Again by applying Lemma~\ref{lem:Pf-Dodgeson} to $A = Y(I \sqcup K)$, 
we obtain the desired formula (\ref{eq:main2}) for $Z(I,J,K)$.

Next we consider the case where $m \ge 4$.
In this case we apply Lemma~\ref{lem:Pf-Dodgeson} to the matrix $A = \tilde{Y} = Z(I,J,K)$ 
with $(i,j,k,l) = (l+1,l+2,l+3,l+4)$.
We see that
\begin{gather*}
\tilde{Y}^{l+\alpha,l+\beta}
 =
Z(I, J \setminus ( j_\alpha, j_\beta ), ( j_\alpha, j_\beta ) \sqcup K)
\quad\text{if $(\alpha, \beta) = (1,2)$, $(2,3)$ or $(3,4)$},
\\
\tilde{Y}^{l+1,l+2,l+3,l+4}
 =
Z(I, J \setminus ( j_1, j_2, j_3, j_4 ), ( j_1, j_2, j_3, j_4 ) \sqcup K),
\end{gather*}
and
\begin{itemize}
\item
$\tilde{Y}^{l+1,l+3}$ is equal to the matrix 
obtained from $-Z(I, J \setminus ( j_1, j_3 ), ( j_1, j_3 ) \sqcup K)$ 
by multiplying the $(l+1)$st row/column (corresponding to $j_2$) by $-1$,
\item
$\tilde{Y}^{l+2,l+4}$ is equal to the matrix 
obtained from $-Z(I, J \setminus ( j_2, j_4 ), ( j_2, j_4 ) \sqcup K)$ 
by multiplying the $(l+2)$nd row/column (corresponding to $j_3$) by $-1$,
\item
$\tilde{Y}^{l+1,l+4}$ is equal to the matrix 
obtained from $Z(I, J \setminus ( j_1, j_4 ), ( j_1, j_4 ) \sqcup K)$ 
by multiplying the $(l+1)$st and $(l+2)$nd rows/columns (corresponding to $j_2$ and $j_3$) by $-1$.
\end{itemize}
Hence, by using the induction hypothesis on $m$, we have
\begin{align*}
&
Z(I,J,K)
\\
 &\quad
 =
\frac{ 1 }
     { \Pf Y( I \sqcup (j_1, j_2, j_3, j_4) \sqcup K ) }
\left\{
 \begin{array}{l}
  \Pf Y( I \sqcup (j_1, j_2) \sqcup K )
  \cdot
  \Pf Y( I \sqcup (j_3, j_4) \sqcup K )
 \\
  -
  \Pf Y( I \sqcup (j_1, j_3) \sqcup K )
  \cdot
  \Pf Y( I \sqcup (j_2, j_4) \sqcup K )
 \\
  +
  \Pf Y( I \sqcup (j_1, j_4) \sqcup K )
  \cdot
  \Pf Y( I \sqcup (j_2, j_3) \sqcup K )
 \end{array}
\right\}
\\
 &\quad\quad
 \times
 \bigl( \Pf Y( J \sqcup K ) \bigr)^{(l+m)/2-1}.
\end{align*}
Again by using Lemma~\ref{lem:Pf-Dodgeson} to $A = Y( I \sqcup (j_1, j_2, j_3, j_4) \sqcup K )$ 
with $(i,j,k,l) = (l+1,l+2,l+3,l+4)$, we obtain the desired formula (\ref{eq:main2}) for $Z(I,J,K)$.

Now it remains to show the cases where
$$
(l,m) = (1,1), \quad (2,2), \quad (1,3), \quad (3,1) \quad\text{and}\quad (3,3).
$$
The case $(l,m) = (1,1)$ is trivial.
The case $(l,m) = (2,2)$ follows from (\ref{eq:special22}) and the cases $(l,m) = (1,3)$ and $(3,1)$ follow 
from (\ref{eq:special13}) in Lemma~\ref{lem:special}.
Hence it is enough to consider the case $(l,m) = (3,3)$.

Let $l = m = 3$, $I = ( a, b, c )$, $J = ( d, e, f )$, 
and write $\langle i_1, \cdots, i_m \rangle = \Pf Y( ( i_1, \cdots, i_m ) \sqcup K)$.
Then we have
$$
Z(I,J,K)
=
\begin{pmatrix}
 0 & \langle a,b,d,e,f \rangle & \langle a,c,d,e,f \rangle 
& \langle a,e,f \rangle & \langle a,d,f \rangle & \langle a,d,e \rangle \\
   & 0                         & \langle b,c,d,e,f \rangle 
& \langle b,e,f \rangle & \langle b,d,f \rangle & \langle b,d,e \rangle \\
   &                           & 0 
& \langle c,e,f \rangle & \langle c,d,f \rangle & \langle c,d,e \rangle \\
   &                           & 
& 0                     & \langle f \rangle     & \langle e \rangle \\
   &                           & 
&                       & 0                     & \langle d \rangle \\
   &                           & 
&                       &                       & 0
\end{pmatrix}.
$$
We apply Lemma~\ref{lem:Pf-Dodgeson} to this matrix, say $W$.
It follows from (\ref{eq:special13}) that
$$
\Pf W^{1,2}
 =
\langle c \rangle \cdot \langle d,e,f \rangle,
\quad
\Pf W^{1,3}
 =
\langle b \rangle \cdot \langle d,e,f \rangle,
\quad
\Pf W^{2,3}
 =
\langle a \rangle \cdot \langle d,e,f \rangle.
$$
And it follows from (\ref{eq:special22}) that
$$
\Pf W^{3,4}
 =
\langle a,b,d \rangle \cdot \langle d,e,f \rangle,
\quad
\Pf W^{2,4}
 =
\langle a,c,d \rangle \cdot \langle d,e,f \rangle,
\quad
\Pf W^{1,4}
 =
\langle b,c,d \rangle \cdot \langle d,e,f \rangle.
$$
Hence, again by using (\ref{eq:special13}), we see that
\begin{align*}
&
\Pf W^{1,2} \cdot \Pf W^{3,4} - \Pf W^{1,3} \cdot \Pf W^{2,4} + \Pf W^{1,4} \cdot \Pf W^{2,3}
\\
&\quad
 =
\langle d,e,f \rangle^2
 \cdot
\left(
 \langle a \rangle \cdot \langle b,c,d \rangle
 - \langle b \rangle \cdot \langle a,c,d \rangle
 + \langle c \rangle \cdot \langle a,b,d \rangle
\right)
\\
&\quad
 =
\langle d,e,f \rangle^2 \cdot \langle a,b,c \rangle \cdot \langle d \rangle.
\end{align*}
Since $\Pf W^{1,2,3,4} = \langle d \rangle$, we conclude that
$$
\Pf W
 =
\langle a,b,c \rangle \cdot \langle d,e,f \rangle^2
 =
\Pf Y( ( a,b,c ) \sqcup K) \cdot \bigl( \Pf Y( ( d,e,f ) \sqcup K) \bigr)^2.
$$
This completes the proof of Theorem~\ref{thm:main2}.
\end{demo}

%%%%%%%%%%
% Proof of generalized Sylvester identity for determinants
%%%%%%%%%%

\section{%
Proof of Theorem~\ref{thm:main1}
}

In this section we derive Theorem~\ref{thm:main1} from Theorem~\ref{thm:main1}.
Some determinant formulas can be deduced from Pfaffian formulas by using the following fundamental relation 
between determinants and Pfaffians.

\begin{lemma}
(See e.g. \cite{IsO}.)
\label{lem:det=Pf}
Suppose that $m+n$ is even.
If $M$ is an $m \times n$ matrix, then we have
\begin{equation}
\label{eq:det=Pf1}
\Pf \begin{pmatrix}
 O & M \\ -\trans M & O
\end{pmatrix}
 =
\begin{cases}
 (-1)^{\binom{m}{2}} \det M &\text{if $m=n$,} \\
 0 &\text{if $m \neq n$.}
\end{cases}
\end{equation}
More generally, if $N$ is an $n \times n$ skew-symmetric matrix and $M$ is an $m \times n$ matrix, then we have
\begin{equation}
\label{eq:det=Pf2}
\Pf \begin{pmatrix}
 N & M \\ -\trans M & O
\end{pmatrix}
 =
\begin{cases}
 (-1)^{\binom{m}{2}} \det M &\text{if $m=n$,} \\
 0 &\text{if $m > n$.}
\end{cases}
\end{equation}
\end{lemma}

\begin{demo}{Proof of Theorem~\ref{thm:main1}}
We put
\begin{alignat*}{3}
A &= X \begin{pmatrix} I \\ I' \end{pmatrix},
&\quad
B &= X \begin{pmatrix} I \\ J' \end{pmatrix},
&\quad
P &= X \begin{pmatrix} I \\ K' \end{pmatrix},
\\
C &= X \begin{pmatrix} J \\ I' \end{pmatrix},
&\quad
D &= X \begin{pmatrix} J \\ J' \end{pmatrix},
&\quad
Q &= X \begin{pmatrix} J \\ K' \end{pmatrix},
\\
R &= X \begin{pmatrix} K \\ I' \end{pmatrix},
&\quad
S &= X \begin{pmatrix} K \\ J' \end{pmatrix},
&\quad
T &= X \begin{pmatrix} K \\ K' \end{pmatrix},
\end{alignat*}
that is,
$$
X = \bordermatrix{
   & I & J & K \cr
I' & A & B & P \cr
J' & C & D & Q \cr
K' & R & S & T
},
$$
and consider the following $2(p+q+r) \times 2(p+q+r)$ skew-symmetric matrix $Y$ with rows/columns indexed by 
$I \sqcup I' \sqcup J \sqcup J' \sqcup K \sqcup K'$:
$$
Y = \bordermatrix{
   & I & I' & J & J' & K & K' \cr
I  & O & A & O & B & O & P \cr
I' & -\trans A & O & -\trans C & O & -\trans R & O \cr
J  &  O & C & O & D & O & Q \cr
J' & -\trans B & O & -\trans D & O & -\trans S & O \cr
K  & 0 & R & O & S & O & T \cr
K' & -\trans P & O & -\trans Q & 0 & -\trans T & O
}.
$$
We apply Theorem~\ref{thm:main2} to the sequences 
$\tilde{I} = I \sqcup I'$, $\tilde{J} = J \sqcup J'$ and $\tilde{K} = K \sqcup K'$.

We compute the entries of the matrix $\tilde{Y}$ given in Theorem~\ref{thm:main2}.
By permuting row/columns and then by using (\ref{eq:det=Pf1}), 
we have for $1 \le \alpha, \beta \le p$
\begin{align*}
&
\Pf Y( (i_\alpha, i_\beta) \sqcup \tilde{J} \sqcup \tilde{K} )
\\
 &\quad
=
(-1)^{qr}
\Pf \begin{pmatrix}
O_{q+r+2}
 &
X \begin{pmatrix} (i_\alpha, i_\beta) \sqcup J \sqcup K \\ J' \sqcup K' \end{pmatrix}
\\
- \trans \left( X \begin{pmatrix} (i_\alpha, i_\beta) \sqcup J \sqcup K \\ J' \sqcup K' \end{pmatrix} \right)
 &
O_{q+r}
\end{pmatrix}
\\
 &\quad
=
0,
\\
&
\Pf Y( (i_\alpha, i'_\beta) \sqcup \tilde{J} \sqcup \tilde{K})
\\
 &\quad
=
(-1)^{qr + (q+r)}
\Pf \begin{pmatrix}
 O_{q+r+1}
&
 X \begin{pmatrix} (i_\alpha) \sqcup J \sqcup K \\ (i'_\beta) \sqcup J' \sqcup K' \end{pmatrix}
\\
 - \trans \left(
 X \begin{pmatrix} (i_\alpha) \sqcup J \sqcup K \\ (i'_\beta) \sqcup J' \sqcup K' \end{pmatrix}
 \right)
&
 O_{q+r+1}
\end{pmatrix}
\\
 &\quad
=
(-1)^{qr + (q+r)} \cdot (-1)^{ \binom{q+r+1}{2} }
\det X \begin{pmatrix} (i_\alpha) \sqcup J \sqcup K \\ (i'_\beta) \sqcup J' \sqcup K' \end{pmatrix},
\\
&
\Pf Y( (i'_\alpha, i'_\beta) \sqcup \tilde{J} \sqcup \tilde{K})
\\
 &\quad
=
(-1)^{r^2 + q(q+r)}
\Pf \begin{pmatrix}
 O_{q+r+2}
&
 -\trans \left(
  X \begin{pmatrix} J \sqcup K \\ (i'_\alpha, i'_\beta) \sqcup J' \sqcup K' \end{pmatrix}
 \right)
\\
 X \begin{pmatrix} J \sqcup K \\ (i'_\alpha, i'_\beta) \sqcup J' \sqcup K' \end{pmatrix}
&
 O_{q+r}
\end{pmatrix}
\\
 &\quad
=
0,
\end{align*}
where $O_n$ denotes the $n \times n$ zero matrix.
Similarly we have, for $1 \le \alpha \le p$ and $1 \le \beta \le q$,
\begin{align*}
\Pf Y( (i_\alpha) \sqcup (\tilde{J} \setminus (j_\beta)) \sqcup \tilde{K})
 &=
(-1)^{ qr + \binom{q+r}{2} } \det X \begin{pmatrix} (i_\alpha) \sqcup (J \setminus (j_\beta)) \sqcup K \\ J' \sqcup K' \end{pmatrix},
\\
\Pf Y( (i_\alpha) \sqcup (\tilde{J} \setminus (j'_\beta)) \sqcup \tilde{K})
 &=
0,
\\
\Pf Y( (i'_\alpha) \sqcup (\tilde{J} \setminus (j_\beta)) \sqcup \tilde{K})
 &=
0,
\\
\Pf Y( (i'_\alpha) \sqcup (\tilde{J} \setminus (j'_\beta)) \sqcup \tilde{K})
 &=
(-1)^{q + qr} \cdot (-1)^{\binom{q+r}{2}}
\det X \begin{pmatrix} J \sqcup K \\ (i'_\alpha) \sqcup (J' \setminus (j'_\beta)) \sqcup K' \end{pmatrix},
\end{align*}
and, for $1 \le \alpha, \beta \le q$, 
\begin{align*}
\Pf Y( (\tilde{J} \setminus ( j_\alpha, j_\beta )) \sqcup \tilde{K})
 &=
0,
\\
\Pf Y( (\tilde{J} \setminus (j_\alpha, j'_\beta)) \sqcup \tilde{K})
 &=
(-1)^{(q-1)r} \cdot (-1)^{\binom{q+r-1}{2}}
\det X \begin{pmatrix} (J \setminus (j_\alpha)) \sqcup K \\ (J' \setminus (j'_\beta)) \sqcup K' \end{pmatrix},
\\
\Pf Y( (\tilde{J} \setminus ( j'_\alpha, j'_\beta )) \sqcup \tilde{K})
 &=
0.
\end{align*}
If we put
\begin{align*}
U
 &= 
\left(
 \det X \begin{pmatrix} (i_\alpha) \sqcup J \sqcup K \\ (i'_\beta) \sqcup J' \sqcup K' \end{pmatrix}
\right)_{1 \le \alpha, \beta \le p},
\\
V
 &=
\left(
 \det X \begin{pmatrix} (i_\alpha) \sqcup (J \setminus (j_\beta)) \sqcup K \\ J' \sqcup K' \end{pmatrix}
\right)_{1 \le \alpha \le p, 1 \le \beta \le q},
\\
V'
 &=
\left(
 \det X \begin{pmatrix} J \sqcup K \\  (i'_\alpha) \sqcup (J' \setminus (j'_\beta)) \sqcup K' \end{pmatrix}
\right)_{1 \le \alpha \le p, 1 \le \beta \le q},
\\
W
 &=
\left(
 \det X \begin{pmatrix} (J \setminus (j_\alpha)) \sqcup K \\ (J' \setminus (j'_\beta)) \sqcup K' \end{pmatrix}
\right)_{1 \le \alpha, \beta \le q},
\end{align*}
then the matrix $\tilde{Y}$ in Theorem~\ref{thm:main2} is given by
$$
\tilde{Y}
 =
\begin{pmatrix}
O              & \ep U                  & \ep V                    & O \\
- \ep \trans U & O                      & O                        & (-1)^q \ep V' \\
- \ep \trans V & O                      & O                        & (-1)^{q-1} \ep W \\
O              & - (-1)^q \ep \trans V' & -(-1)^{q-1} \ep \trans W & O
\end{pmatrix},
$$
where $\ep = (-1)^{qr+\binom{q+r}{2}}$.
By permuting rows/columns and using Lemma~\ref{lem:det=Pf}, we have
\begin{align*}
\Pf \tilde{Y}
 &=
(-1)^{q(p+q)} \cdot (-1)^{\binom{p+q}{2}}
\det \begin{pmatrix}
\ep U                 & \ep V                    \\
-(-1)^q \ep \trans V' & -(-1)^{q-1} \ep \trans W
\end{pmatrix}
\\
 &=
(-1)^{q(p+q) + \binom{p+q}{2} + (p+q)(qr+ \binom{q+r}{2})} \cdot (-1)^{q+q^2}
\det \begin{pmatrix}
U         & V          \\
\trans V' & - \trans W
\end{pmatrix}.
\end{align*}
Since we have
\begin{align*}
\Pf Y(\tilde{I} \sqcup \tilde{K})
 &=
(-1)^{pr+\binom{p+r}{2}} \det X \begin{pmatrix}I \sqcup K \\ I' \sqcup K' \end{pmatrix},
\\
\Pf Y(\tilde{J} \sqcup \tilde{K})
 &=
(-1)^{qr+ \binom{q+r}{2}} \det X \begin{pmatrix} J \sqcup K \\ J' \sqcup K' \end{pmatrix},
\end{align*}
we obtain
\begin{multline*}
(-1)^{q(p+q) + \binom{p+q}{2} + (p+q) ( qr + \binom{q+r}{2} ) + q + q^2}
\det \tilde{X}
 \\
 =
(-1)^{pr + \binom{p+r}{2} + (p+q-1) ( qr + \binom{q+r}{2})}
\det X \begin{pmatrix} I \sqcup K \\ I' \sqcup K' \end{pmatrix}
\cdot
\left( \det X \begin{pmatrix} J \sqcup K \\ J' \sqcup K' \end{pmatrix} \right)^{p+q-1}.
\end{multline*}
Now by noticing the congruence
$$
pq + \binom{p+q}{2} 
+ qr + \binom{q+r}{2}
+ rp + \binom{r+p}{2}
 \equiv
0
\ \bmod 2,
$$
we complete the proof of Theorem~\ref{thm:main1}.
\end{demo}

We can derive the following Bazin formula from Theorem~\ref{thm:main1}

\begin{corollary}
(Bazin \cite{B})
Let $n$ and $p$ be positive integers such that $p \le n$.
If $Z$ is a matrix with columns indexed by $( 1, 2, \dots, n )$, 
and $I = (i_1, \dots, i_p)$, $J = (j_1, \dots, j_p)$, and $K$ are sequences of row indices 
of length $p$, $p$, and $n-p$ respectively, 
then we have
\begin{multline}
\label{eq:Bazin}
\det \left(
 \det Z \begin{pmatrix} (i_\alpha) \sqcup (J \setminus (j_\beta)) \sqcup K \\ (1, 2, \dots, n) \end{pmatrix}
\right)_{1 \le \alpha, \beta \le p}
\\
 =
\det Z \begin{pmatrix} I \sqcup K \\ (1, 2, \dots, n) \end{pmatrix}
\cdot
\left( \det Z \begin{pmatrix} J \sqcup K \\ (1, 2, \dots, n) \end{pmatrix} \right)^{p-1}.
\end{multline}
\end{corollary}

\begin{demo}{Proof}
We apply Theorem~\ref{thm:main1} to the matrix $X = \left( x_{i,j} \right)$ 
with columns indexed by $(\overline{1}, \dots, \overline{p}, \allowbreak 1, \dots, p, p+1, \dots, n)$ given by
$$
\begin{cases}
 x_{i,\overline{j}} = z_{i,j} &\text{if $1 \le j \le p$,} \\
 x_{i,j} = z_{i,j} &\text{if $1 \le j \le n$.}
\end{cases}
$$
Let $I' = (\overline{1}, \dots, \overline{p})$, $J' = ( 1, \dots, p)$, $K' = (p+1, \dots, n)$.
Then we have
\begin{gather*}
\det X \begin{pmatrix} (i_\alpha) \sqcup J \sqcup K \\ (\overline{i}) \sqcup J' \sqcup K' \end{pmatrix}
 =
0
\quad(1 \le \alpha \le p, 1 \le i \le p),
\\
\det X \begin{pmatrix} (i_\alpha) \sqcup (J \setminus (j_\beta)) \sqcup K \\ J' \sqcup K' \end{pmatrix}
 =
\det Z \begin{pmatrix} (i_\alpha) \sqcup (J \setminus (j_\beta)) \sqcup K \\ (1, 2, \dots, n) \end{pmatrix}
\quad(1 \le \alpha, \beta \le p),
\\
\det X \begin{pmatrix} J \sqcup K \\ (\overline{i}) \sqcup (J' \setminus (j)) \sqcup K' \end{pmatrix}
 =
\delta_{i,j} (-1)^{i-1} \det Z \begin{pmatrix} J \sqcup K \\ (1, 2, \dots, n) \end{pmatrix}
\quad(1 \le i, j \le p),
\\
\det X \begin{pmatrix} I \sqcup K \\ I' \sqcup K' \end{pmatrix}
 =
\det Z \begin{pmatrix} I \sqcup K \\ (1, 2, \dots, n) \end{pmatrix},
\\
\det X \begin{pmatrix} J \sqcup K \\ J' \sqcup K' \end{pmatrix}
 =
\det Z \begin{pmatrix} J \sqcup K \\ (1, 2, \dots, n) \end{pmatrix}.
\end{gather*}
Hence we have
\begin{multline*}
\det \begin{pmatrix}
 O_p
&
 \left(
  \det Z \begin{pmatrix} (i_\alpha) \sqcup (J \setminus (j_\beta)) \sqcup K \\ (1, 2, \dots, n) \end{pmatrix}
 \right)_{1 \le \alpha, \beta \le p}
\\
 \left(
  \delta_{i,j} (-1)^{i-1} \det Z \begin{pmatrix} J \sqcup K \\ (1, 2, \dots, n) \end{pmatrix}
 \right)_{1 \le i, j \le p}
&
 *
\end{pmatrix}
\\
 =
(-1)^p
\det Z \begin{pmatrix} I \sqcup K \\ (1, 2, \dots, n) \end{pmatrix}
\cdot
\left(
 \det Z \begin{pmatrix} J \sqcup K \\ (1, 2, \dots, n) \end{pmatrix}
\right)^{2p-1}.
\end{multline*}
By permuting columns on the determinant of the left hand side, we see that the left hand side is equal to
\begin{multline*}
(-1)^{p^2}
\det \left(
 \det Z \begin{pmatrix} (i_\alpha) \sqcup (J \setminus (j_\beta)) \sqcup K \\ (1, 2, \dots, n) \end{pmatrix}
\right)_{1 \le \alpha, \beta \le p}
\\
\times
(-1)^{\sum_{i=1}^p (i-1)} 
\left(
 \det Z \begin{pmatrix} J \sqcup K \\ (1, 2, \dots, n) \end{pmatrix}
\right)^p.
\end{multline*}
By cancelling the common factor, we obtain the Bazin formula (\ref{eq:Bazin}).
\end{demo}

%%%%%%%%%%
% Applications to symmetric functions
%%%%%%%%%%

\section{%
Skew generalizations of Giambelli identity
}

In this section, we use the generalized Sylvester formula (Theorem~\ref{thm:main1}) 
to obtain a skew generalization of the Giambelli identity to skew Schur functions.

A partition of a nonnegative integer $n$ is a weakly decreasing sequence 
$\lambda = (\lambda_1, \lambda_2, \lambda_3, \dots)$ 
of nonnegative integers such that $|\lambda| = \sum_{i \ge 1} \lambda_i = n$.
The length, denoted by $l(\lambda)$, of a partition $\lambda$ is the number of nonzero entries of $\lambda$.
We sometimes write $\lambda = (\lambda_1, \dots, \lambda_{l(\lambda)})$ by omitting the $0$s at the end.
We identify a partition $\lambda$ of $n$ with its Young diagram, 
which is a left-justified array of $n$ cells with $\lambda_i$ cells in the $i$th row.
We denote by $\emptyset$ the empty partition $(0,0,\dots)$ of $0$.
Given a partition $\lambda$, we put
$$
p(\lambda) = \# \{ i : \lambda_i \ge i \},
\quad
\alpha_i = \lambda_i - i,
\quad
\beta_i = \lambda'_i - i
\quad(1 \le i \le p(\lambda)),
$$
where $\lambda'_i$ is the number of cells in the $i$th column of the Young diagram of $\lambda$.
Then we write $\lambda = (\alpha_1, \dots, \alpha_{p(\lambda)} | \beta_1, \dots, \beta_{p(\lambda)})$ 
and call it the Frobenius notation of $\lambda$.

Let $s_\lambda$ be the Schur function corresponding to a partition $\lambda$, 
and $s_{\lambda/\mu}$ the skew Schur function associated to a pair of partitions $\lambda$ and $\mu$.
Note that $s_\lambda = s_{\lambda/\emptyset}$ and $s_{\mu/\mu} = 1$ 
and that $s_{\lambda/\mu} = 0$ unless $\lambda_i \ge \mu_i$ for $i \ge 1$.
(Refer to \cite[Chapter~I]{Macdonald} for details on Schur functions.)
Giambelli \cite{G} gave a formula which expresses any Schur function as a determinant of Schur functions 
of hook shapes $(a|b)$, and Lascoux--Pragacz \cite{LP} gave a generalization of the Giambelli identity 
to skew Schur functions.

\begin{prop}
\label{prop:schur_det}
\begin{enumerate}
\item[(1)]
(Giambelli \cite{G})
For a partition $\lambda = (\alpha_1, \dots, \alpha_p | \beta_1, \dots, \beta_p)$, we have
\begin{equation}
\label{eq:Giambelli}
s_\lambda = \det \left( s_{(\alpha_i | \beta_j)} \right)_{1 \le i, j \le p}.
\end{equation}
\item[(2)]
(Lascoux--Pragacz \cite{LP})
If $\lambda = (\alpha_1, \cdots, \alpha_p|\beta_1, \cdots, \beta_p)$ 
and $\mu = (\gamma_1, \cdots, \gamma_q|\delta_1, \cdots, \delta_q)$ in the Frobenius notation, then we have
\begin{equation}
\label{eq:LP}
s_{\lambda/\mu}
 =
(-1)^q
\det \begin{pmatrix}
 \Bigl(
  s_{(\alpha_i|\beta_j)}
 \Bigr)_{1 \le i, j \le p}
&
 \Bigl(
  h_{\alpha_i - \gamma_j}
 \Bigr)_{\substack{1 \le i \le p \\ 1 \le j \le q}}
\\
 \Bigl( 
  e_{\beta_j - \delta_i}
 \Bigr)_{\substack{1 \le i \le q \\ 1 \le j \le p}}
&
 O
\end{pmatrix},
\end{equation}
where $h_k$ (resp. $e_k$) is the $k$th complete (resp. elementary) symmetric function 
and $h_k = e_k = 0$ for $k < 0$.
\end{enumerate}
\end{prop}

We apply Theorem~\ref{thm:main1} to give another generalization of 
the Giambelli identity to skew Schur functions.
In order to state the identity, we extend the definition of skew Schur functions as follows:

\begin{definition}
Given two nonnegative integer sequences $\alpha = (\alpha_1, \dots, \alpha_p)$ 
and $\beta = (\beta_1, \dots, \beta_p)$ of the same length $p$ 
and a partition $\mu$, we define 
$s_{(\alpha|\beta)/\mu}$ by putting
$$
s_{(\alpha|\beta)/\mu}
 =
(-1)^q
\det \begin{pmatrix}
 \Bigl(
  s_{(\alpha_i|\beta_j)}
 \Bigr)_{1 \le i, j \le p}
&
 \Bigl( 
  h_{\alpha_i - \gamma_j}
 \Bigr)_{\substack{1 \le i \le p \\ 1 \le j \le q}}
\\
 \Bigl(
  e_{\beta_j - \delta_i}
 \bigr)_{\substack{1 \le i \le q \\ 1 \le j \le p}}
&
 O
\end{pmatrix},
$$
where $\mu = (\gamma_1, \cdots, \gamma_q|\delta_1, \cdots, \delta_q)$ in the Frobenius notation.
\end{definition}

If the entries of $\alpha$ (or $\beta$) are not distinct, then $s_{(\alpha|\beta)/\mu} = 0$.
Otherwise, if $\sigma$ and $\tau \in S_p$ are permutations such that
$\alpha_{\sigma(1)} > \cdots > \alpha_{\sigma(p)}$ and $\beta_{\tau(1)} > \cdots > \beta_{\tau(p)}$, 
then by (\ref{eq:LP}) we have $s_{(\alpha|\beta)/\mu} = \sgn(\sigma \tau) s_{\lambda/\mu}$, 
where $\lambda$ is a partition given by the Frobenius notation 
$\lambda = (\alpha_{\sigma(1)}, \cdots, \alpha_{\sigma(p)} | \beta_{\tau(1)}, \cdots, \beta_{\tau(p)})$.
Then we have the following skew generalization of the Giambelli identity.

\begin{theorem}
\label{thm:skewGiambelli}
If two partitions $\lambda$ and $\mu$ are represented as 
$\lambda = (\alpha_1, \cdots, \alpha_p|\beta_1, \cdots, \beta_p)$ and 
$\mu = (\gamma_1, \cdots, \gamma_q|\delta_1, \cdots, \delta_q)$ in the Frobenius notation, 
then we have
\begin{multline}
\label{eq:skewGiambelli}
s_{\lambda/\mu}
\\
 =
(-1)^q
\det
\begin{pmatrix}
 \Bigl(
  s_{(\alpha_i,\gamma_1, \cdots, \gamma_q|\beta_j,\delta_1, \cdots, \delta_q)/\mu}
 \Bigr)_{1 \le i, j \le p}
&
 \bigl(
  s_{(\alpha_i, \gamma_1, \cdots, \widehat{\gamma_j}, \cdots, \gamma_q|\delta_1, \cdots, \delta_q)/\mu}
 \Bigr)_{\substack{1 \le i \le p \\ 1 \le j \le q}}
\\
 \Bigl(
  s_{(\gamma_1, \cdots, \gamma_q|\beta_j, \delta_1, \cdots, \widehat{\delta_i}, \cdots, \delta_q)/\mu}
 \Bigr)_{\substack{1 \le i \le q \\ 1 \le j \le p}}
&
 O
\end{pmatrix},
\end{multline}
where the symbol $\hat{a}$ means removing $a$ from the sequence.
\end{theorem}

If $\mu = \emptyset$, then (\ref{eq:skewGiambelli}) reduces to the Giambelli identity (\ref{eq:Giambelli}).
Note that the nonzero entries of the determinant (\ref{eq:skewGiambelli}) are of the form 
$\pm s_{\nu/\mu}$ with $\nu/\mu$ a border strip, i.e., a connected skew Young diagram 
containing no $2 \times 2$ block of cells.

\begin{demo}{Proof}
Apply Theorem~\ref{thm:main1} to the matrix
$$
X = \begin{pmatrix}
 \Bigl( s_{(\alpha_i|\beta_j)} \Bigr)_{1 \le i, j \le p}
&
 \Bigl( s_{(\alpha_i|\delta_j)} \Bigr)_{\substack{1 \le i \le p \\ 1 \le j \le q}}
&
 \Bigl( h_{\alpha_i-\gamma_j} \Bigr)_{\substack{1 \le i \le p \\ 1 \le j \le q}}
\\
 \Bigl( s_{(\gamma_i|\beta_j)} \Bigr)_{\substack{1 \le i \le q \\ 1 \le j \le p}}
&
 \Bigl( s_{(\gamma_i|\delta_j)} \Bigr)_{1 \le i, j \le q}
&
 \Bigl( h_{\gamma_i-\gamma_j} \Bigr)_{1 \le i, j \le q}
\\
 \Bigl( e_{\beta_j-\delta_i} \Bigr)_{\substack{1 \le i \le q \\ 1 \le j \le p}}
&
 \Bigl( e_{\delta_j-\delta_i} \Bigr)_{1 \le i, j \le q}
&
 O
\end{pmatrix}
$$
with $I = I' = ( 1, \cdots, p )$, $J = J' = ( p+1, \cdots, p+q )$ and $K = K' = ( p+q+1, \cdots, p+2q )$.
Then by definition we have
\begin{gather*}
\det X \begin{pmatrix} ( i ) \sqcup J \sqcup K \\ ( j ) \sqcup J \sqcup K \end{pmatrix}
 =
(-1)^q s_{(\alpha_i,\gamma_1, \cdots, \gamma_q|\beta_j,\delta_1, \cdots, \delta_q)/\mu},
\\
\det X \begin{pmatrix} ( i ) \sqcup (J \setminus ( j )) \sqcup K \\ J \sqcup K \end{pmatrix}
 =
(-1)^q s_{(\alpha_i, \gamma_1, \cdots, \hat{\gamma_j}, \cdots, \gamma_q|\delta_1, \cdots, \delta_q)/\mu},
\\
\det X \begin{pmatrix} J \sqcup K \\ ( i ) \sqcup (J \setminus ( j )) \sqcup K \end{pmatrix}
 =
(-1)^q s_{(\gamma_1, \cdots, \gamma_q| \beta_i, \delta_1, \cdots, \hat{\delta_j}, \cdots, \delta_q)/\mu},
\\
\det X \begin{pmatrix} I \sqcup K \\ I \sqcup K \end{pmatrix}
 =
(-1)^q s_{(\alpha|\beta)/\mu}
 =
(-1)^q s_{\lambda/\mu},
\\
\det X \begin{pmatrix} J \sqcup K \\ J \sqcup K \end{pmatrix}
 =
(-1)^q s_{(\gamma|\delta)/\mu}
 =
(-1)^q s_{\mu/\mu}
 =
(-1)^q.
\end{gather*}
Since $X \begin{pmatrix} (J \setminus ( i )) \sqcup K \\ (J \setminus ( j )) \sqcup K \end{pmatrix}$ 
is a $(2s-1) \times (2s-1)$ matrix whose bottom-right block is the $s \times s$ zero matrix, 
we see that
$$
\det X \begin{pmatrix} (J \setminus ( i )) \sqcup K \\ (J \setminus ( j )) \sqcup K \end{pmatrix}
 =
0.
$$
Hence, by applying Theorem~\ref{thm:main1}, we obtain the desired identity.
\end{demo}

\begin{remark}
\label{rem:KP}
As is shown in \cite{NOS}, Theorem~\ref{thm:skewGiambelli} is obtained by using 
the Giambelli-type determinant formula for the expansion coefficients 
of the $\tau$-function $\tau(x)$ of the KP hierarchy (\cite[Theorem~1.1]{NOS}) 
and the fact that $\tau(x) = \sum_\lambda s_{\lambda/\mu}(y) s_\lambda(x)$ is a solution of 
the KP-hierarchy.
\end{remark}

%%%%%%%%%%
% Applications to symmetric functions
%%%%%%%%%%

\section{%
Skew generalizations of Schur identity
}

In this section we use Theorem~\ref{thm:main2} to obtain a skew generalization of the Schur Pfaffian identity 
for Schur's $Q$-functions.

A partition $\lambda$ is called strict if $\lambda_1 > \lambda_2 > \dots > \lambda_{l(\lambda)}$.
Let $Q_\lambda$ be the Schur $Q$-function corresponding to a strict partition $\lambda$, 
and $Q_{\lambda/\mu}$ the skew $Q$-function associated with a pair of strict partitions $\lambda$ and $\mu$.
Note that $Q_{\lambda/\emptyset} = Q_\lambda$, $Q_{\mu/\mu} = 1$ and 
$Q_{\lambda/\mu} = 0$ unless $\lambda_i \ge \mu_i$ for all $i$.
See \cite[Chapter~III, Section~8]{Macdonald} for details on Schur $Q$-functions.

Schur \cite{Schur} defined the Schur $Q$-function corresponding to any strict partition 
as a Pfaffian of Schur $Q$-functions corresponding to strict partitions with at most two rows,
and J\'ozefiak--Pragacz (\cite{JP}) gave a generalization of Schur's identity to 
skew $Q$-functions.

\begin{prop}
\label{prop:Schur}
For sequences $\alpha = (\alpha_1, \dots, \alpha_l)$ and $\beta = (\beta_1, \dots, \beta_m)$ 
of nonnegative integers, let $S^\alpha_\beta$ and $T^\alpha_\beta$ be 
$l \times m$ matrices defined by
$$
S^\alpha_\beta
 =
\Bigl( Q_{(\alpha_i,\beta_j)} \Bigr)_{1 \le i \le l, 1 \le j \le m},
\quad
T^\alpha_\beta
 =
\Bigl( Q_{(\alpha_i - \beta_{m+1-j})} \Bigr)_{1 \le i \le l, 1 \le j \le m},
$$
where we use the convention
$$
Q_{(a,b)} = - Q_{(b,a)},
\quad
Q_{(a,0)} = - Q_{(0,a)} = Q_{(a)},
\quad
Q_{(0,0)} = 0
$$
for positive integers $a$ and $b$, and $Q_{(a)} = 0$ for $a < 0$.
Then we have the following Pfaffian identities:
\begin{enumerate}
\item[(1)]
(Schur \cite{Schur})
If $\lambda$ is a strict partition, then we have
\begin{equation}
\label{eq:Schur}
Q_\lambda
 = 
\begin{cases}
 \Pf S^\lambda_\lambda &\text{if $l(\lambda)$ is even,} \\
 \Pf S^{\lambda^0}_{\lambda^0} &\text{if $l(\lambda)$ is odd,}
\end{cases}
\end{equation}
where $\lambda = (\lambda_1, \dots, \lambda_{l(\lambda)})$ and 
$\lambda^0 = (\lambda_1, \dots, \lambda_{l(\lambda)}, 0)$.
\item[(2)]
(J\'ozefiak--Pragacz \cite{JP})
For two strict partitions $\lambda$ and $\mu$, we have
\begin{equation}
\label{eq:JP}
Q_{\lambda/\mu}
 =
\begin{cases}
\Pf \begin{pmatrix}
 S^\lambda_\lambda & T^\lambda_\mu \\
 -\trans T^\lambda_\mu & O
\end{pmatrix}
 &\text{if $l(\lambda) + l(\mu)$ is even,} \\
\Pf \begin{pmatrix}
 S^\lambda_\lambda & T^\lambda_{\mu^0} \\
 -\trans T^\lambda_{\mu^0} & O
\end{pmatrix}
 &\text{if $l(\lambda) + l(\mu)$ is odd.}
\end{cases}
\end{equation}
\end{enumerate}
\end{prop}

We find another generalization of the Schur identity (\ref{eq:Schur}) to skew $Q$-functions.
To state the identity, we extend the definition of skew $Q$-functions as follows:

\begin{definition}
Given a nonnegative integer sequence $\alpha$ of length $l$ and a strict partition $\mu$, 
we define $Q_{\alpha/\mu}$ by putting
\begin{align*}
Q_{\alpha/\mu}
 &=
\begin{cases}
\Pf \begin{pmatrix}
 S^\alpha_\alpha & T^\alpha_\mu \\
 -\trans T^\alpha_\mu & O
\end{pmatrix}
 &\text{if $l + l(\mu)$ is even,} \\
\Pf \begin{pmatrix}
 S^\alpha_\alpha & T^\alpha_{\mu^0} \\
 -\trans T^\alpha_{\mu^0} & O
\end{pmatrix}
 &\text{if $l + l(\mu)$ is odd.}
\end{cases}
\end{align*}
\end{definition}

We note that, if $l+l(\mu)$ is odd, then we have
$$
\begin{pmatrix}
 S^\alpha_\alpha & T^\alpha_{\mu^0} \\
 -\trans T^\alpha_{\mu^0} & O
\end{pmatrix}
 =
\begin{pmatrix}
 S^{\alpha^0}_{\alpha^0} & T^{\alpha^0}_\mu \\
 -\trans T^{\alpha^0}_\mu & O
\end{pmatrix},
$$
where $\alpha^0 = (\alpha_1, \dots, \alpha_l, 0)$.
If the entries of $\alpha$ are not distinct, then $Q_{\alpha/\mu} = 0$.
Otherwise, if $\sigma \in S_l$ is a permutation such that $\alpha_{\sigma(1)} > \cdots > \alpha_{\sigma(l)}$, 
then $Q_{\alpha/\mu} = \sgn(\sigma) Q_{\lambda/\mu}$ with $\lambda = (\alpha_{\sigma(1)}, \cdots, \alpha_{\sigma(l)})$.
Then we have the following skew-generalization of Schur identity (\ref{eq:Schur}).

\begin{theorem}
\label{thm:skewSchur}
For strict partitions $\lambda$ and $\mu$, we have
\begin{equation}
\label{eq:skewSchur}
Q_{\lambda/\mu}
 =
\Pf \begin{pmatrix}
 \Bigl(
  Q_{(\lambda_i,\lambda_j,\mu_1, \cdots, \mu_m)/\mu}
 \Bigr)_{1 \le i, j \le l}
&
 \Bigl(
  Q_{(\lambda_i, \mu_1, \cdots, \hat{\mu_j}, \cdots, \mu_s)/\mu}
 \Bigr)_{\substack{1 \le i \le l \\ 1 \le j \le s}}
\\
 - \trans \Bigl(
  Q_{(\lambda_i, \mu_1, \cdots, \hat{\mu_j}, \cdots, \mu_s)/\mu}
 \Bigr)_{\substack{1 \le i \le l \\ 1 \le j \le s}}
&
 O
\end{pmatrix},
\end{equation}
where $l = l(\lambda)$, $m = l(\mu)$ and $s = m$ or $m+1$ according to whether $l+m$ is even or odd.
\end{theorem}

If $\mu = \emptyset$, then (\ref{eq:skewSchur}) reduces to the Schur identity (\ref{eq:Schur}).

\begin{demo}{Proof}
We put $l = l(\lambda)$ and $m = l(\mu)$.

First we consider the case where $l \equiv m \bmod 2$.
In this case we apply Theorem~\ref{thm:main2} to the $(l+2m) \times (l+2m)$ skew-symmetric matrix
$$
Y
 =
\begin{pmatrix}
 S^\lambda_\lambda & S^\lambda_\mu & T^\lambda_\mu \\
 -\trans S^\lambda_\mu & S^\mu_\mu & T^\mu_\mu \\
 -\trans T^\lambda_\mu & -\trans T^\mu_\mu & O
\end{pmatrix},
$$
with $I = \{ 1, \cdots, l \}$, $J = \{ l+1, \cdots, l+m \}$ and $K = \{ l+m+1, \cdots, l+2m \}$.
Then by definition we have
\begin{gather*}
\Pf Y( ( i, j ) \sqcup J \sqcup K )
 =
Q_{(\lambda_i, \lambda_j, \mu_1, \cdots, \mu_m)/\mu},
\\
\Pf Y( ( i ) \sqcup (J \setminus ( l+j )) \sqcup K )
 =
Q_{(\lambda_i, \mu_1, \cdots, \hat{\mu_j}, \cdots, \mu_m)/\mu}.
\end{gather*}
Since $Y( (J \setminus ( l+i, l+j )) \sqcup K)$ is a $(2m-2) \times (2m-2)$ skew-symmetric matrix 
whose bottom-right block is the $m \times m$ zero matrix, we have 
$\Pf Y( (J \setminus ( l+i, l+j )) \sqcup K) = 0$ by (\ref{eq:det=Pf2}).
Also we have
$$
\Pf Y( I \sqcup K) = Q_{\lambda/\mu},
\quad
\Pf Y( J \sqcup K) = Q_{\mu/\mu} = 1.
$$
Hence, by applying (\ref{eq:main2}), we obtain (\ref{eq:skewSchur}).

Next we consider the cases where $l \not\equiv m \bmod 2$.
In this case , we apply Theorem~\ref{thm:main2} to the $(l+1+2m) \times (l+1+2m)$ skew-symmetric matrix
$$
Y
 =
\begin{pmatrix}
 S^{\lambda^0}_{\lambda^0} & S^{\lambda^0}_\mu & T^{\lambda^0}_\mu \\
 -\trans S^{\lambda^0}_\mu & S^\mu_\mu & T^\mu_\mu \\
 -\trans T^{\lambda^0}_\mu & -\trans T^\mu_\mu & O
\end{pmatrix},
$$
with $I = \{ 1, \cdots, l+1 \}$, $J = \{ l+2, \cdots, l+m+1 \}$ and $K = \{ l+m+2, \cdots, l+2m+2 \}$ 
and $\lambda_{l+1} = 0$.
Then by definition we have
\begin{gather*}
\Pf Y( ( i, j ) \sqcup J \sqcup K )
 =
Q_{(\lambda_i,\lambda_j,\mu_1, \cdots, \mu_m)/\mu}
\quad\text{if $1 \le i < j \le l$,}
\\
\Pf Y( (i) \sqcup (J \setminus (l+1+j)) \sqcup K )
 =
Q_{(\lambda_i, \mu_1, \dots, \hat{\mu_j}, \dots, \mu_m)/\mu}
\quad\text{if $1 \le i \le l$ and $1 \le j \le m$.}
\end{gather*}
If $1 \le i \le l$, then by moving the $2$nd column/row to the $(m+1)$st column/row we see that
\begin{align*}
\Pf Y( ( i, l+1 ) \sqcup J \sqcup K )
 &=
(-1)^m
\Pf \begin{pmatrix}
 S^{(\lambda_i, \mu_1, \dots, \mu_m)}_{(\lambda_i,\mu_1, \dots, \mu_m)}
&
 T^{(\lambda_i, \mu_1, \dots, \mu_m)}_{(\mu_1, \dots, \mu_m, 0)}
\\
 -\trans T^{(\lambda_i, \mu_1, \dots, \mu_m)}_{(\mu_1, \dots \mu_m,0)}
&
 O
\end{pmatrix}
\\
 &=
(-1)^m Q_{(\lambda_i,\mu_1, \cdots, \mu_m)/\mu}.
\end{align*}
If $1 \le j \le m$, then by moving the $1$st row/column to the $m$th row/column 
and by noting $Q_{(\lambda_{l+1}-\mu_k)} = 0$ for $1 \le k \le m$, we have
$$
\Pf Y( ( l+1 ) \sqcup (J \setminus ( l+1+ j ) \sqcup K )
 =
(-1)^{m-1}
\Pf \begin{pmatrix}
 S^{(\mu_1, \dots, \hat{\mu_j}, \dots, \mu_m)}_{(\mu_1, \dots, \hat{\mu_j}, \dots, \mu_m)}
&
 T^{(\mu_1, \dots, \hat{\mu_j}, \dots, \mu_m)}_{(\mu_1, \dots, \mu_m,0)}
\\
 -\trans T^{(\mu_1, \dots, \hat{\mu_j}, \dots, \mu_m)}_{(\mu_1, \dots, \mu_m,0)}
&
 O_{m+1}
\end{pmatrix}.
$$
Hence by using (\ref{eq:det=Pf2}), we have 
$\Pf Y( ( l+1 ) \sqcup (J \setminus ( l+1+ j ) \sqcup K ) = 0$.
Also we have
$$
\Pf Y( I \sqcup K) = Q_{\lambda/\mu},
\quad
\Pf Y( J \sqcup K) = Q_{\mu/\mu} = 1.
$$
Therefore by applying (\ref{eq:main2}), we obtain
$$
Q_{\lambda/\mu} = \Pf \tilde{Y},
$$
where the entries of the skew-symmetric matrix $\tilde{Y} = \bigl( \tilde{y}_{i,j} \bigr)_{1 \le i, j \le l+m+1}$ 
are given by
$$
\begin{cases}
\tilde{y}_{i,j} = Q_{(\lambda_i,\lambda_j,\mu_1, \cdots, \mu_m)/\mu} &\text{if $1 \le i, j \le l$,} \\
\tilde{y}_{i,l+1} = (-1)^m Q_{(\lambda_i,\mu_1, \cdots, \mu_m)/\mu} &\text{if $1 \le i \le l$,} \\
\tilde{y}_{i,l+1+j} = Q_{(\lambda_i, \mu_1, \cdots, \hat{\mu_j}, \cdots, \mu_m)/\mu}
 &\text{if $1 \le i \le l$ and $1 \le j \le m$,} \\
\tilde{y}_{i,j} = 0 &\text{if $l+1 \le i, j \le l+m+1$.}
\end{cases}
$$
By pulling out the common factor $(-1)^m$ from the $(l+1)$st row/column 
and moving the $(l+1)$st row/column to the last row/column, we see that
\begin{multline*}
\Pf \tilde{Y}
 =
(-1)^m \cdot (-1)^m
\\
\times
\Pf \begin{pmatrix}
 \Bigl(
  Q_{(\lambda_i,\lambda_j,\mu_1, \cdots, \mu_m)/\mu}
 \Bigr)_{1 \le i, j \le l}
&
 \Bigl(
  Q_{(\lambda_i, \mu_1, \cdots, \hat{\mu_j}, \cdots, \mu_{m+1})/\mu}
 \Bigr)_{\substack{1 \le i \le l \\ 1 \le j \le m+1}}
\\
 - \trans \Bigl(
  Q_{(\lambda_i, \mu_1, \cdots, \hat{\mu_j}, \cdots, \mu_{m+1})/\mu}
 \Bigr)_{\substack{1 \le i \le l \\ 1 \le j \le m+1}}
&
 O
\end{pmatrix}.
\end{multline*}
This completes the proof of Theorem~\ref{thm:skewSchur}.
\end{demo}

\begin{remark}
\label{rem:BKP}
(This remark is due to A.~Nakayashiki.) 
Theorem~\ref{thm:skewSchur} can be also obtained from the theory of the BKP hierarchy 
by using the same idea as in \cite{NOS} (see Remark~\ref{rem:KP}).
We fix a strict partition $\mu$ and consider a formal power series of the form
$$
\tau(x)
 =
Q_\mu \left( \frac{x}{2} \right)
 +
\sum_\lambda \xi_\lambda Q_\lambda \left( \frac{x}{2} \right),
$$
where $\lambda$ runs over strict partitions with $|\lambda| > |\mu|$ 
and $x = (x_1, x_3, x_5, \dots)$ is the so-called Sato variables, 
i.e., $x_i = p_i/i$ in the symmetric function language.
Then Shigyo \cite[Theorem~3]{Shigyo} proves that $\tau(x)$ is a solution of the BKP hierarchy if and only if
the coefficients $\xi_\lambda$ satisfy the following Pfaffian formulas:
$$
\xi_\lambda
 =
\Pf \begin{pmatrix}
 \Bigl(
  \xi_{(\lambda_i,\lambda_j,\mu_1, \cdots, \mu_m)}
 \Bigr)_{1 \le i, j \le l}
&
 \Bigl(
  \xi_{(\lambda_i, \mu_1, \cdots, \hat{\mu_j}, \cdots, \mu_s)}
 \Bigr)_{\substack{1 \le i \le l \\ 1 \le j \le s}}
\\
 - \trans \Bigl(
  \xi_{(\lambda_i, \mu_1, \cdots, \hat{\mu_j}, \cdots, \mu_s)}
 \Bigr)_{\substack{1 \le i \le l \\ 1 \le j \le s}}
&
 O
\end{pmatrix},
$$
where $l = l(\lambda)$, $m=l(\mu)$, and $s = m$ or $m+1$ according to whether $l+m$ is even or odd, 
and we use the convention that $\xi_{(\alpha_1, \dots, \alpha_p)}$ is alternating 
in $\alpha_1, \dots, \alpha_p$.
By using the fact that $Q_\mu(x/2)$ is a solution of the BKP hierarchy \cite{You}, we can show that 
$$
\tau(x)
 = 
\sum_\lambda Q_{\lambda/\mu}(y) Q_\lambda \left( \frac{x}{2} \right)
 =
Q_\mu \left( \frac{x}{2} \right)
\exp \left( \sum_n n x_n y_n \right),
$$
where $n$ runs over all positive odd integers and $y = (y_1, y_3, y_5, \dots)$ is another set of variables, 
is a solution of the BKP hierarchy.
Applying Shigyo's formula to this special solution $\tau(x)$, we obtain the generalized Schur identity 
(\ref{eq:skewSchur}).
\end{remark}

\section*{%
Acknowledgment
}

This work is motivated by questions raised by Prof.~A.~Nakayashiki and Dr.~Y.~Shigyo.
The author thanks them for valuable discussions and comments.

\end{document}